\documentstyle[12pt]{amsart}
\def\NZQ{\Bbb}               % the font for N,Z,Q,R,C

\def\QQ{{\NZQ Q}}
\def\ZZ{{\NZQ Z}}
\def\RR{{\NZQ R}}

\def\PP{{\NZQ P}}
\def\ab{{\bold a}}

\def\eb{{\bold e}}

\def\xb{{\bold x}}

\def\opn#1#2{\def#1{\operatorname{#2}}} % to make operators
\opn\chara{char}
\opn\length{\ell}
\opn\pd{pd}
\opn\rk{rk}
\opn\projdim{proj\,dim}
\opn\injdim{inj\,dim}
\opn\rank{rank}
\opn\depth{depth}
\opn\grade{grade}
\opn\height{height}
\opn\embdim{emb\,dim}
\opn\codim{codim}

\opn\Tr{Tr}
\opn\bigrank{big\,rank}
\opn\superheight{superheight}\opn\lcm{lcm}
\opn\trdeg{tr\,deg}%
\opn\reg{reg}
\opn\lreg{lreg}
\opn\skel{skel}
\opn\com{com}
\opn\div{div}
\opn\Div{Div}
\opn\cl{cl}
\opn\Cl{Cl}
\opn\Spec{Spec}
\opn\Supp{Supp}
\opn\supp{supp}
\opn\Sing{Sing}
\opn\Ass{Ass}
\opn\Ann{Ann}
\opn\Rad{Rad}
\opn\Soc{Soc}
\opn\Ker{Ker}
\opn\Coker{Coker}
\opn\Im{Im}
\opn\Hom{Hom}
\opn\Tor{Tor}
\opn\Ext{Ext}
\opn\End{End}
\opn\Aut{Aut}
\opn\id{id}

\opn\nat{nat}
\opn\pff{pf}%   \pf exists already
\opn\Pf{Pf}
\opn\GL{GL}
\opn\SL{SL}
\opn\mod{mod}
\opn\ord{ord}
\opn\aff{aff}
\opn\con{conv}
\opn\relint{relint}
\opn\st{st}
\opn\lk{lk}
\opn\cn{cn}
\opn\core{core}
\opn\vol{vol}
\opn\link{link}
\opn\star{star}
\opn\gr{gr}

\def\pot#1#2{#1[\kern-0.28ex[#2]\kern-0.28ex]}

\opn\dirlim{\underrightarrow{\lim}}
\opn\inivlim{\underleftarrow{\lim}}
\let\Union=\bigcup
\let\Sect=\bigcap

\let\to=\rightarrow

\def\Implies{\ifmmode\Longrightarrow \else
     \unskip${}\Longrightarrow{}$\ignorespaces\fi}
\def\implies{\ifmmode\Rightarrow \else
     \unskip${}\Rightarrow{}$\ignorespaces\fi}
\def\iff{\ifmmode\Longleftrightarrow \else
     \unskip${}\Longleftrightarrow{}$\ignorespaces\fi}

\let\:=\colon
\newtheorem{Theorem}{Theorem}[section]

\newtheorem{Corollary}[Theorem]{Corollary}

\newtheorem{Example}[Theorem]{Example}

\let\epsilon\varepsilon
\let\phi=\varphi
\let\kappa=\varkappa
\def\qed{\ifhmode\textqed\fi
   \ifmmode\ifinner\quad\qedsymbol\else\dispqed\fi\fi}
\def\textqed{\unskip\nobreak\penalty50
    \hskip2em\hbox{}\nobreak\hfil\qedsymbol
    \parfillskip=0pt \finalhyphendemerits=0}
\def\dispqed{\rlap{\qquad\qedsymbol}}

\def\HH{{\cal H}}

\def\PP{{\cal P}}
\def\QQ{{\cal Q}}

\opn\inii{in}
\opn\inim{inm}
\opn\rate{rate}
\opn\Mon{Mon}
\opn\lex{lex}
\opn\rev{rev}
\opn\red{red}
\opn\max{max}
\opn\min{min}
\opn\initial{in}
\opn\incom{incom}
\opn\com{com}
\opn\Ker{Ker}
\opn\Krull{Krull}
\begin{document}
\title{Special simplices and Gorenstein toric rings}
\author{Hidefumi Ohsugi and Takayuki Hibi}
\date{}
\maketitle
\begin{abstract}
Christos Athanasiadis \cite{Christos} studies an effective
technique to show that
Gorenstein sequences 
coming from
% the Ehrhart polynomial of 
compressed polytopes are 
% $(0,1)$-polytopes are
unimodal. 
In the present paper we will use such the technique to find
a rich class of Gorenstein toric rings 
with unimodal $h$-vectors arising from finite graphs.
\end{abstract}

\section*{Introduction}
Let $\PP \subset \RR^N$ be an integral convex polytope, 
i.e., a convex polytope 
each of whose vertices
has integer coordinates.
Let 
$K[\xb, \xb^{-1}, t] = 
K[x_1, x_1^{-1}, \ldots, x_N, x_N^{-1}, t]$
denote the Laurent polynomial ring in $(N + 1)$ variables over 
a field $K$.
The {\em toric ring} of $\PP$ is
the subalgebra $K[\PP]$
of $K[\xb, \xb^{-1}, t]$
generated by those Laurent polynomials
$\xb^{\ab}t = x_1^{a_1} \cdots x_N^{a_N} t$
such that $\ab = (a_1, \ldots, a_N)$
is a vertex of $\PP$.
% We say that $K[\PP]$ is the {\em toric ring} of $\PP$.
We will regard $K[\PP]$ as a homogeneous algebra
\cite[p. 147]{BruHer}
by setting each $\deg \xb^{\ab}t = 1$
and write
$F(K[\PP], \lambda)$ for its Hilbert series. 
One has 
$F(K[\PP], \lambda) = 
(h_0 + h_1 \lambda + \cdots + h_s \lambda^s) / (1 - \lambda)^{d+1}$,
where each $h_i \in \ZZ$ with $h_s \neq 0$ and where $d$ is the
dimension of $\PP$.  The sequence $(h_0, h_1, \ldots, h_s)$
is said to be the $h$-{\em vector} of $K[\PP]$.
% Recall that $\PP$ is {\em normal}
If the toric ring $K[\PP]$
is normal,   
% If $\PP$ is normal, 
then $K[\PP]$ is Cohen--Macaulay.
If $K[\PP]$ is Cohen--Macaulay, then
the $h$-vector of $K[\PP]$ is nonnegative,
i.e., each $h_i \geq 0$.  Moreover, if $K[\PP]$ is Gorenstein,
then the $h$-vector of $K[\PP]$ is symmetric, i.e.,
$h_i = h_{s-i}$ for all $i$.

% Let $V(\PP)$ denote the set of vertices of $\PP$
% and $K[\yb] = K[\{ y_\ab \}_{\ab \in V(\PP)}]$ denote
% the polynomial ring in $|V(\PP)|$ variables over $K$ with 
% each $\deg y_\ab = 1$.  The {\em toric ideal} of $K[\PP]$ 
% is the kernel $I_{\PP}$ of the surjective ring homomorphism
% $\pi : K[\yb] \to K[\PP]$ defined by setting
% $\pi(y_\ab) = \xb^{\ab}t$ for all $\ab \in V(\PP)$.
% We say that $\PP$ is {\em compressed} 
% if each of the reverse lexicographic initial ideals of $I_\PP$ 
% is generated by squarefree monomials.

An outstanding conjecture (which is still open) is
that the $h$-vector of a Gorenstein toric ring is unimodal, 
i.e., $h_0 \leq h_1 \leq \cdots \leq h_{[s/2]}$.
One of the established techniques to show that
the $h$-vector $(h_0, h_1, \ldots, h_s)$
of a Gorenstein toric ring $K[\PP]$ is unimodal is 
to find a simplicial convex polytope
\cite{Stanley}
% \cite[pp. 000]{StanleyGreenBook}
of dimension $s - 1$ 
whose $h$-vector 
coincides with
$(h_0, h_1, \ldots, h_s)$.
On the other hand, however, 
given a Gorenstein toric ring $K[\PP]$,
it seems difficult
to find such a simplicial convex polytope.

Christos Athanasiadis \cite{Christos} introduces the concept 
of a special simplex of a convex polytope.  
Let $\PP \subset \RR^N$ be a convex polytope.
A $(q-1)$-simplex $\Sigma$ each of whose vertices is
a vertex of $\PP$ is said to be a {\em special simplex} 
in $\PP$ if each facet of $\PP$ contains exactly $q - 1$ of
the vertices of $\Sigma$.
Recall that
an integral convex polytope $\PP \subset \RR^N$
is {\em compressed} \cite[p. 337]{Sta}
(and \cite{OhHicompressed})
if all
``pulling triangulations'' 
of $\PP$ are unimodular.  
The toric ring $K[\PP]$ of a compressed polytope $\PP$ is normal. 
It turns out \cite[Theorem 3.5]{Christos}
that
if $\PP$ is compressed and if
there is a special simplex in $\PP$, then
% $\PP$ possesses a special simplex, then 
the $h$-vector of $K[\PP]$ is equal to the $h$-vector of 
a simplicial convex polytope.

In the present paper we will use \cite[Theorem 3.5]{Christos}
to study the $h$-vector of the toric ring of the edge polytope of 
a finite graph satisfying the odd cycle condition 
as well as that of 
the stable polytope of a perfect graph.

\newpage

\section{Two polytopes arising from finite graphs}
Let $G$ be a finite graph on the vertex set 
$[n] = \{ 1, 2, \ldots, n \}$
having no loops and no multiple edges,
and $E(G)$ the edge set of $G$.
We associate each subset $W \subset [n]$ with 
the $(0, 1)$-vector
$\rho(W) = \sum_{j \in W} \eb_j \in \RR^n$.
Here $\eb_j$ is the $j$-th unit coordinate vector 
in $\RR^n$.  Thus in particular $\rho(\emptyset)$ 
is the origin of $\RR^n$.
A subset $W \subset [n]$ is called {\em stable}
(resp. a {\em clique})
if $\{i, j \} \not\in E(G)$ 
(resp. $\{i, j \} \in E(G)$)
for all $i, j \in W$
with $i \neq j$.
Note that the empty set as well as each single-element 
susbset of $[n]$ is both stable and a clique.
Let $S(G)$ denote the set of stable sets of $G$.

We now introduce two convex polytopes arising 
from a finite graph $G$ on $[n]$.  
First, the {\em edge polytope} \cite{OhHinormal} of $G$ 
is the $(0, 1)$-polytope $\PP_G \subset \RR^n$ 
which is the convex hull of 
$\{ \rho(e) \, : \, e \in E(G) \}$.
Second, the {\em stable polytope} \cite{Chvatal}
of $G$ is the $(0, 1)$-polytope $\QQ_G \subset \RR^n$ 
which is the convex hull of 
$\{ \rho(W) \, : \, W \in S(G) \}$.

% The {\em clique complex} of $G$ on 
% is the simplicial complex $\Delta(G)$ on $[n]$
% whose faces are the cliques of $G$. 

\begin{Example}
\label{posetpolytope}
{\em
Let $P$ be a finite poset on $[n]$ 
and $\com(P)$ its comparability graph.
Thus $\com(P)$ is the finite graph on $[n]$ 
such that $\{ i , j \}$ with $i \neq j$
is an edge of $\com(P)$ if and only if
$i$ and $j$ are comparable in $P$.
Then the stable polytope of $\com(P)$ coincides with
the chain polytope \cite{twoposetpolytopes} of $P$.
}
\end{Example}

% In \cite{OhHinormal} (\cite{SVV})
The problem
when the toric ring $K[\PP_G]$ is normal
and the problem 
when the edge polytope $\PP_G$ possesses a unimodular 
covering \cite[p. 420]{OhHinormal} are studied
in \cite{OhHinormal} (and \cite{SVV}).

\begin{Theorem}[\cite{OhHinormal}]
\label{normal}
Given a finite connected graph $G$,
the following conditions are equivalent:
\begin{enumerate}
\item[(i)]
% \cite[Corollary 2.3]{OhHinormal} and
% \cite[Theorem 00.00]{SVV}
The toric ring $K[\PP_G]$ is normal;
\item[(ii)] 
The edge polytope $\PP_G$ possesses a unimodular covering;
\item[(iii)]
$G$ satisfies the odd cycle condition, i.e., 
if each of $C$ and $C'$  
is an odd cycle (a cycle of odd length)
of $G$ and if $C$ and $C'$ possess no common vertex, 
then there exists an edge $\{ i, j \}$ 
of $G$ such that $i$ 
is a vertex of $C$
and $j$ is a vertex of $C'$. 
\end{enumerate}
Thus in particular the edge polytope of
a finite connected bipartite graph 
possesses a unimodular covering 
and its toric ring is normal.
\end{Theorem}

A {\em chromatic number} of a finite graph $G$ on $[n]$ 
is the smallest integer $\ell > 0$ for which there is a map 
$\varphi : [n] \to [\ell]$
with the property that $\varphi(i) \neq \varphi(j)$
if $\{i, j \} \in E(G)$. 
A finite graph $G$ is called {\em perfect} 
if, for all induced subgraphs $H$ of $G$ including $G$ itself,
the chromatic number of $H$ is equal to the maximal 
cardinality of cliques contained in $H$.
The comparability graph of a finite 
partially ordered set is perfect (\cite{Berge}).

% A convex polytope is called integral
% if each of its vertices has integer coordinates.
% A {\em compressed polytope} \cite[p. 337]{Sta}
% is an integral convex polytope all of whose
% ``pulling triangulations'' are unimodular.  
% We refer the reader to \cite{OhHicompressed} 
% for the detailed information about compressed polytopes.

The facets of the edge polytope $\PP_G$ of a finite connected
graph $G$ is completely determined 
(\cite[Theorem 1.7]{OhHinormal}).
On the other hand, the facets of the stable polytope $\QQ_G$ 
is completely determined when $G$ is a perfect graph
(\cite[Theorem 3.1]{Chvatal}). 
% The stable polytope of a perfect graph is compressed
% (\cite[Example 00.00]{OhHicompressed}).

\section{Gorenstein toric rings}
When $K[\PP_G]$ (resp. $K[\QQ_G]$) is normal, 
it follows easily that $K[\PP_G]$ (resp. $K[\QQ_G]$) 
coincides with the {\em Ehrhart ring} 
\cite[p. 97]{Hibi} of $\PP_G$ (resp. $\QQ_G$). 
On the other hand, 
since the equations of the facets of
$\PP_G$ (resp. $\QQ_G$) are known, 
when $K[\PP_G]$ (resp. $K[\QQ_G]$) is normal,
by using the criterion 
\cite[Corollary (1.2)]{DeNegri--Hibi}
one can determines the finite graphs $G$ for which 
the toric ring 
$K[\PP_G]$ (resp. $K[\QQ_G]$)
is Gorenstein. 

Let $G$ be a finite connected graph on $[n]$.  
Given a subset $V \neq \emptyset$ of $[n]$, write $G_V$
for the induced subgraph of $G$ on $V$.  
We say that $G$ is {\em $2$-connected} if 
$G$ together with $G_{[n] \setminus \{ i \}}$
for all $i \in [n]$ is connected.
If $i \in [n]$, then $N(G;i)$ stands for the set of
vertices $j$ with $\{ i, j \} \in E(G)$.
If $T \subset [n]$, then 
$N(G;T) = \Union_{i \in T} N(G;i)$.
The {\em bipartite graph induced by} 
a stable set $T \neq \emptyset$ of $G$
is the bipartite graph on the vertex set $T \Union N(G;T)$ 
consisting of those edges $\{ i, j \}$ of $G$ 
with $i \in T$ and $j \in N(G;T)$.

Recall that a {\em matching} of 
$G$ is a set of edges $\{ e_1, \ldots, e_m \}$
such that $e_i \Sect e_j = \emptyset$ for all $i \neq j$.
A matching $\{ e_1, \ldots, e_m \}$ of $G$ is called
{\em perfect} if $\Union_{i=1}^{m} e_i = [n]$.
In particular $n$ is even and $m = n / 2$
if $G$ possesses a perfect matching
$\{ e_1, \ldots, e_m \}$.
It follows that $G$ possesses a perfect matching
if and only if the monomial $x_1 x_2 \cdots x_n$ belongs to 
the toric ring $K[\PP_G]$. 

\begin{Theorem}
\label{Gorenstein}
{\em (a)} 
Let $G$ be a finite connected graph on $[n]$
satisfying the odd cycle condition  
and suppose that every connected component of 
$G_{[n] \setminus \{ i \}}$
possesses at least one odd cycle for all $i \in [n]$.
Then the toric ring $K[\PP_G]$ of the edge polytope $\PP_G$
of $G$ is Gorenstein if and only if
(i) $G$ possesses a perfect matching,
(ii)
one has
$|N(G;T)| = |T| + 1$ for each  
stable set $T$ of $G$ such that
the bipartite graph induced by $T$ is connected
with
$T \Union N(G;T) \neq [n]$ and that
every connected component of 
$G_{[n] \setminus (T \Union N(G;T))}$ 
has at least one odd cycle
and (iii)
one has
$|T| = n/2 - 1$ for each 
stable set $T$ of $G$ such that
the bipartite graph induced by $T$ is connected
with
$T \Union N(G;T) = [n]$.

% \smallskip

{\em (a')} 
Let $G$ be a bipartite graph on 
$[n] = V_1 \Union V_2$
and suppose that $G$ is $2$-connected.  
Then the toric ring $K[\PP_G]$ of the edge polytope $\PP_G$
of $G$ is Gorenstein if and only if 
(i) $G$ possesses a perfect matching and (ii)
one has
$|N(G;T)| = |T| + 1$ 
for every subset 
$T \subset V_1$ 
such that
$G_{T \Union N(G;T)}$ is connected and that 
$G_{[n] \setminus (T \Union N(G;T))}$
is a connected graph with at least one edge. 

% \smallskip

{\em (b)} 
The toric ring $K[\QQ_G]$ of a stable polytope $\QQ_G$ of 
a perfect graph $G$ is Gorenstein if and only if 
all maximal cliques have the same cardinality. 
\end{Theorem}

\begin{pf}
% We give a proof of (a) by using the criterion 
% \cite[Corollary (1.2)]{DeNegri--Hibi} 
% and omit a proof of each of (a') and (b) 
% which will be done similarly.
%
% \smallskip
%
(a)  The edge polytope
$\PP_G \subset \RR^n$ 
lies on the hyperplane $\HH$
defined by the equation
$z_1 + \cdots + z_n = 2$.
Let $\pi : \RR^{n-1} \to \HH$ denote the affine map 
defined by setting
$\pi(z_1, \ldots, z_{n-1}) 
= (z_1, \ldots, z_{n-1}, 2 - (z_1 + \cdots + z_{n-1}))$.
Then $\pi$ is an affine isomorphism with
$\pi(\ZZ^{n-1}) = \HH \Sect \ZZ^n$.
Hence $\pi^{-1}(\PP_G) \subset \RR^{n-1}$ is an integral
convex polytope with $\dim \pi^{-1}(\PP_G) = n - 1$ and 
the toric ring $K[\pi^{-1}(\PP_G)]$ is isomorphic to
$K[\PP_G]$ as homogeneous algebras over $K$.

Let $\delta$ denote the smallest integer for which
the interior of
% 
% $\delta \PP_G$
%
$\delta (\pi^{-1}(\PP_G))$
contains at least one integer point
$(a_1, \ldots, a_{n-1})$.
Since every connected component of 
$G_{[n] \setminus \{ i \}}$
possesses at least one odd cycle, 
it follows from
\cite[Theorem 1.7 (a)]{OhHinormal}
that the hyperplane defined by 
the equation $z_i = 0$  
is a supporting hyperplane 
which defines a facet of $\PP_G$
% $\pi^{-1}(\PP_G)$ 
for each $1 \leq i \leq n$.
Thus the hyperplane defined by 
the equation $z_i = 0$  
is a supporting hyperplane 
which defines a facet of 
$\pi^{-1}(\PP_G)$ 
for each $1 \leq i < n$.
Thus by using \cite[Corollary (1.2)]{DeNegri--Hibi}
one has each $a_i = 1$
if $K[\pi^{-1}(\PP_G)]$ is Gorenstein.
If $\eb_1 + \cdots + \eb_{n-1}$
belongs to the interior of $\delta (\pi^{-1}(\PP_G))$,
then $\eb_1 + \cdots + \eb_{n-1} + q\eb_n$
belongs to the interior of $\delta \PP_G$
for some integer $q > 0$.
Since $K[\PP_G]$ coincides with the Ehrhart ring
of $\PP_G$, it follows that 
there are edges $e_1, \ldots, e_m$ of $G$
with
$\eb_1 + \cdots + \eb_{n-1} + q\eb_n
= \rho(e_1) + \cdots \rho(e_m)$.
Hence $q = 1$ and 
$\eb_1 + \cdots + \eb_n
= \rho(e_1) + \cdots \rho(e_m)$.
Thus $G$ possesses a perfect matching
with $\delta = n / 2$
if $K[\pi^{-1}(\PP_G)]$ is Gorenstein. 

For a while, suppose that 
$G$ possesses a perfect matching
with $\delta = n / 2$. 
Let $\PP^{\flat} \subset \RR^{n-1}$ denote the integral
convex polytope 
$\delta (\pi^{-1}(\PP_G)) - (\eb_1 + \cdots + \eb_{n-1})$.
Then $\PP^{\flat}$ is of standard type, i.e.,
$\dim \PP = n - 1$ and the origin of $\RR^{n-1}$
belongs to the interior of $\PP$.
Then \cite[Corollary (1.2)]{DeNegri--Hibi}
guarantees that the toric ring $K[\pi^{-1}(\PP_G)]$ 
is Gorenstein if and only if 
the dual polytope \cite[p. 631]{DeNegri--Hibi} 
of $\PP^{\flat}$ is integral.

Now, by using \cite[Theorem 1.7 (a)]{OhHinormal}
again,
it turns out that the equations of the supporting 
hyperplanes 
which defines the facets of 
$\PP^{\flat}$ are the followings:  
\begin{itemize}
\item
% [(i)]
$z_i = - 1$ \, for each \, $1 \leq i < n$;
\item
% [(ii)]
$\sum_{i \in [n] \setminus (T \Union N(G;T))} z_i 
+ 2 \sum_{j \in N(G;T)} z_j = |T| - |N(G;T)|$
\, if \, $n \in T$;
\item
% [(iii)]
$2 \sum_{i \in T} z_i 
+ \sum_{j \in [n] \setminus (T \Union N(G;T))} z_j = 
|N(G;T)| - |T|$
\, if \, $n \in N(G;T)$;
\item
% [(iv)]
$\sum_{i \in T} z_i - \sum_{j \in N(G;T)} z_j = 
|N(G;T)| - |T|$
\, if \, $n \not\in T \Union N(G;T)$,
\end{itemize}
where $T \neq \emptyset$ is a stable set of $G$ 
for which the bipartite graph induced by $T$ is connected
and for which either
$T \Union N(G;T) = [n]$
or every connected component of the induced subgraph
$G_{[n] \setminus (T \Union N(G;T))}$ 
% on $[n] \setminus (T \Union N(G;T))$
has at least one odd cycle. 

Hence the dual polytope of $\PP^{\flat}$ is integral 
if and only if ($\alpha$) one has
$|N(G;T)| = |T| + 1$ for each nonempty 
stable set $T$ of $G$ such that
the bipartite graph induced by $T$ is connected
with
$T \Union N(G;T) \neq [n]$ and that
every connected component of 
$G_{[n] \setminus (T \Union N(G;T))}$ 
has at least one odd cycle
and ($\beta$)
one has
$|T| = n/2 - 1$ for each 
stable set $T$ of $G$ such that
the bipartite graph induced by $T$ is connected
with
$T \Union N(G;T) = [n]$.

Consequently, 
when the toric ring $K[\PP_G]$ is Gorenstein,
the conditions (i), (ii) and (iii) are satisfied.
Conversely, suppose that
the conditions (i), (ii) and (iii) are satisfied.
Since the hyperplane defined by the equation $z_i = 0$  
is a supporting hyperplane 
which defines a facet of $\PP_G$
% $\pi^{-1}(\PP_G)$ 
for each $1 \leq i \leq n$,
if $\gamma \PP_G$, where $\gamma > 0$,
contains at least one integer point
$(a_1, \ldots, a_n)$, then
each $a_i > 0$ and $\gamma \geq [(n + 1)/2]$.
It follows from (i), (ii) and (iii) that
$n$ is even and
$\eb_1 + \cdots + \eb_n$ belongs to the interior of
$(n / 2) \PP_G$.
Thus the smallest number $\delta > 0$ for which
$\delta \PP_G$
contains at least one integer point is $\delta = n / 2$
and $\eb_1 + \cdots + \eb_n$ belongs to the interior of
$\delta \PP_G$.
Our discussion done already in the preceding paragraph 
guarantees that the toric ring $K[\PP_G]$ is Gorenstein, 
as desired. 

(a')  In imitation of the preceding proof of (a) 
by using \cite[Theorem 1.7 (b)]{OhHinormal}
instead of \cite[Theorem 1.7 (a)]{OhHinormal},
one can easily give a proof of (a').

(b) The facets of the stable polytope $\QQ_G$ 
is completely determined when $G$ is a perfect graph
(\cite[Theorem 3.1]{Chvatal}). 
In fact, when $G$ is perfect,
the equations of the supporting hyperplanes 
which defines the facets of 
$\QQ_G$ are either $z_i = 0$
for $1 \leq i \leq n$ or
$\sum_{W \subset [n]} z_i = 1$,
where $W$ is a maximal cliques of $G$.
Let $\delta$ denote the smallest integer $\delta > 0$
for which the interior of $\QQ_G$ contains at least
one integer point.  Then $\delta - 1$ coincides with 
the maximal cardinality of cliques of $G$ and
$\eb_1 + \cdots + \eb_n$ belongs to  
the interior of $\delta \QQ_G$.
It follows that 
the dual polytope of the integral polytope
$\delta \QQ_G - (\eb_1 + \cdots + \eb_n)
\subset \RR^n$ of standard type
is integral if and only if
all maximal cliques of $G$ have the cardinality
$\delta - 1$.  Hence  
\cite[Corollary (1.2)]{DeNegri--Hibi}
guarantees that the toric ring $K[\QQ_G]$
is Gorenstein if and only if 
all maximal cliques have the same cardinality. 
\end{pf}

\begin{Example}
\label{graph}
{\em
The toric ring of the edge polytope of each of the finite
connected graphs $G_1$ and $G_2$
drawn below is normal and Gorenstein.
}
\end{Example}

%%%%%%%%%%%%%%% fig %%%%%%%%%%%%%%%%%%%%%%%%%%%%%%%

\bigskip

\bigskip

\bigskip

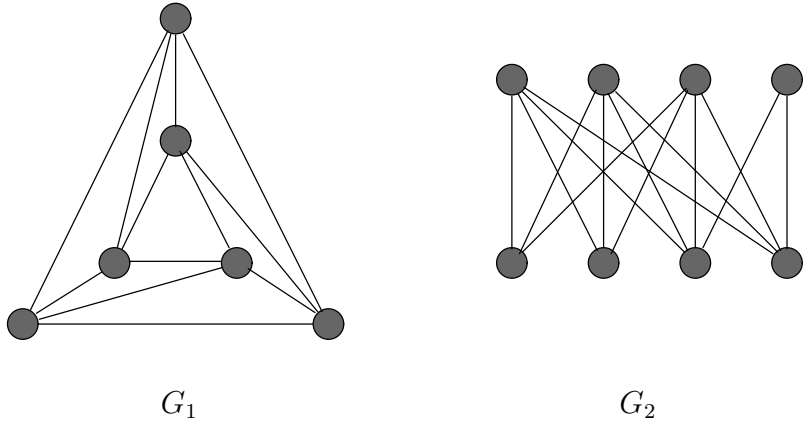
\begin{center}
%WinTpicVersion2.15
\unitlength 0.1in
\begin{picture}(41.60,20.15)(7.00,-29.15)
% CIRCLE 2 0 0 0
% 4 1580 1380 1580 1460 1580 1460 1580 1460
% 
\special{pn 8}%
\special{sh 0.600}%
\special{ar 1580 980 80 80  0.0000000 6.2831853}%
% CIRCLE 2 0 0 0
% 4 1580 2020 1580 2100 1580 2100 1580 2100
% 
\special{pn 8}%
\special{sh 0.600}%
\special{ar 1580 1620 80 80  0.0000000 6.2831853}%
% CIRCLE 2 0 0 0
% 4 1260 2660 1260 2740 1260 2740 1260 2740
% 
\special{pn 8}%
\special{sh 0.600}%
\special{ar 1260 2260 80 80  0.0000000 6.2831853}%
% CIRCLE 2 0 0 0
% 4 1900 2660 1900 2740 1900 2740 1900 2740
% 
\special{pn 8}%
\special{sh 0.600}%
\special{ar 1900 2260 80 80  0.0000000 6.2831853}%
% CIRCLE 2 0 0 0
% 4 2380 2980 2380 3060 2380 3060 2380 3060
% 
\special{pn 8}%
\special{sh 0.600}%
\special{ar 2380 2580 80 80  0.0000000 6.2831853}%
% CIRCLE 2 0 0 0
% 4 3340 1700 3340 1780 3340 1780 3340 1780
% 
\special{pn 8}%
\special{sh 0.600}%
\special{ar 3340 1300 80 80  0.0000000 6.2831853}%
% CIRCLE 2 0 0 0
% 4 780 2980 780 3060 780 3060 780 3060
% 
\special{pn 8}%
\special{sh 0.600}%
\special{ar 780 2580 80 80  0.0000000 6.2831853}%
% CIRCLE 2 0 0 0
% 4 3820 1700 3820 1780 3820 1780 3820 1780
% 
\special{pn 8}%
\special{sh 0.600}%
\special{ar 3820 1300 80 80  0.0000000 6.2831853}%
% CIRCLE 2 0 0 0
% 4 4300 1700 4300 1780 4300 1780 4300 1780
% 
\special{pn 8}%
\special{sh 0.600}%
\special{ar 4300 1300 80 80  0.0000000 6.2831853}%
% CIRCLE 2 0 0 0
% 4 4780 1700 4780 1780 4780 1780 4780 1780
% 
\special{pn 8}%
\special{sh 0.600}%
\special{ar 4780 1300 80 80  0.0000000 6.2831853}%
% CIRCLE 2 0 0 0
% 4 4780 2660 4780 2740 4780 2740 4780 2740
% 
\special{pn 8}%
\special{sh 0.600}%
\special{ar 4780 2260 80 80  0.0000000 6.2831853}%
% CIRCLE 2 0 0 0
% 4 4300 2660 4300 2740 4300 2740 4300 2740
% 
\special{pn 8}%
\special{sh 0.600}%
\special{ar 4300 2260 80 80  0.0000000 6.2831853}%
% CIRCLE 2 0 0 0
% 4 3820 2660 3820 2740 3820 2740 3820 2740
% 
\special{pn 8}%
\special{sh 0.600}%
\special{ar 3820 2260 80 80  0.0000000 6.2831853}%
% CIRCLE 2 0 0 0
% 4 3340 2660 3340 2740 3340 2740 3340 2740
% 
\special{pn 8}%
\special{sh 0.600}%
\special{ar 3340 2260 80 80  0.0000000 6.2831853}%
% LINE 2 0 3 0
% 2 1580 1468 1580 1948
% 
\special{pn 8}%
\special{pa 1580 1068}%
\special{pa 1580 1548}%
\special{fp}%
% LINE 2 0 3 0
% 2 1340 2652 1820 2652
% 
\special{pn 8}%
\special{pa 1340 2252}%
\special{pa 1820 2252}%
\special{fp}%
% LINE 2 0 3 0
% 2 860 2980 2300 2980
% 
\special{pn 8}%
\special{pa 860 2580}%
\special{pa 2300 2580}%
\special{fp}%
% LINE 2 0 3 0
% 2 1540 1444 820 2908
% 
\special{pn 8}%
\special{pa 1540 1044}%
\special{pa 820 2508}%
\special{fp}%
% LINE 2 0 3 0
% 2 1620 1460 2348 2916
% 
\special{pn 8}%
\special{pa 1620 1060}%
\special{pa 2348 2516}%
\special{fp}%
% LINE 2 0 3 0
% 2 3820 1780 3820 2580
% 
\special{pn 8}%
\special{pa 3820 1380}%
\special{pa 3820 2180}%
\special{fp}%
% LINE 2 0 3 0
% 2 4300 1780 4300 2580
% 
\special{pn 8}%
\special{pa 4300 1380}%
\special{pa 4300 2180}%
\special{fp}%
% LINE 2 0 3 0
% 2 4780 1788 4780 2588
% 
\special{pn 8}%
\special{pa 4780 1388}%
\special{pa 4780 2188}%
\special{fp}%
% LINE 2 0 3 0
% 2 1276 2580 1556 1460
% 
\special{pn 8}%
\special{pa 1276 2180}%
\special{pa 1556 1060}%
\special{fp}%
% LINE 2 0 3 0
% 2 1300 2588 1540 2092
% 
\special{pn 8}%
\special{pa 1300 2188}%
\special{pa 1540 1692}%
\special{fp}%
% LINE 2 0 3 0
% 2 1860 2580 1604 2076
% 
\special{pn 8}%
\special{pa 1860 2180}%
\special{pa 1604 1676}%
\special{fp}%
% LINE 2 0 3 0
% 2 1636 2076 2324 2916
% 
\special{pn 8}%
\special{pa 1636 1676}%
\special{pa 2324 2516}%
\special{fp}%
% LINE 2 0 3 0
% 2 852 2924 1196 2708
% 
\special{pn 8}%
\special{pa 852 2524}%
\special{pa 1196 2308}%
\special{fp}%
% LINE 2 0 3 0
% 2 868 2956 1828 2684
% 
\special{pn 8}%
\special{pa 868 2556}%
\special{pa 1828 2284}%
\special{fp}%
% LINE 2 0 3 0
% 2 3340 1788 3340 2588
% 
\special{pn 8}%
\special{pa 3340 1388}%
\special{pa 3340 2188}%
\special{fp}%
% LINE 2 0 3 0
% 2 3372 1772 3788 2596
% 
\special{pn 8}%
\special{pa 3372 1372}%
\special{pa 3788 2196}%
\special{fp}%
% LINE 2 0 3 0
% 2 3844 1772 4260 2596
% 
\special{pn 8}%
\special{pa 3844 1372}%
\special{pa 4260 2196}%
\special{fp}%
% LINE 2 0 3 0
% 2 4332 1772 4748 2596
% 
\special{pn 8}%
\special{pa 4332 1372}%
\special{pa 4748 2196}%
\special{fp}%
% LINE 2 0 3 0
% 2 4740 1772 4340 2580
% 
\special{pn 8}%
\special{pa 4740 1372}%
\special{pa 4340 2180}%
\special{fp}%
% LINE 2 0 3 0
% 2 4260 1756 3860 2612
% 
\special{pn 8}%
\special{pa 4260 1356}%
\special{pa 3860 2212}%
\special{fp}%
% LINE 2 0 3 0
% 2 3404 1740 4708 2612
% 
\special{pn 8}%
\special{pa 3404 1340}%
\special{pa 4708 2212}%
\special{fp}%
% LINE 2 0 3 0
% 2 3396 1764 4244 2604
% 
\special{pn 8}%
\special{pa 3396 1364}%
\special{pa 4244 2204}%
\special{fp}%
% LINE 2 0 3 0
% 2 3884 1764 4732 2604
% 
\special{pn 8}%
\special{pa 3884 1364}%
\special{pa 4732 2204}%
\special{fp}%
% LINE 2 0 3 0
% 2 4252 1748 3396 2604
% 
\special{pn 8}%
\special{pa 4252 1348}%
\special{pa 3396 2204}%
\special{fp}%
% LINE 2 0 3 0
% 2 3780 1756 3380 2612
% 
\special{pn 8}%
\special{pa 3780 1356}%
\special{pa 3380 2212}%
\special{fp}%
% STR 2 0 3 0
% 3 1600 3300 1600 3400 5 0
% $G_1$
\put(16.0000,-30.0000){\makebox(0,0){$G_1$}}%
% STR 2 0 3 0
% 3 4000 3300 4000 3400 5 0
% $G_2$
\put(40.0000,-30.0000){\makebox(0,0){$G_2$}}%
% LINE 2 0 3 0
% 2 1960 2690 2310 2920
% 
\special{pn 8}%
\special{pa 1960 2290}%
\special{pa 2310 2520}%
\special{fp}%
\end{picture}%
\end{center}

\bigskip

\bigskip

\bigskip

%%%%%%%%%%%%%%%%%%%%  fig end %%%%%%%%%%%%%%%%%%%%%%%%%%%%%%%%%%%

\section{Unimodal Gorenstein sequences}
Let $\PP \subset \RR^N$ be a convex polytope.
Recall that 
a $(q-1)$-simplex $\Sigma$ each of whose vertices is
a vertex of $\PP$ is said to be a {\em special simplex}
\cite{Christos} 
in $\PP$ if each facet of $\PP$ contains exactly $q - 1$ of
the vertices of $\Sigma$.

\begin{Theorem}
\label{specialsimplex}
{\em (a)}
Let $G$ be a finite connected graph 
as in Theorem \ref{Gorenstein} (a) or (a')
and suppose that 
the toric ring $K[\PP_G]$ of the edge polytope $\PP_G$
of $G$ is Gorenstein.  
Then there is a special simplex in $\PP_G$.

% \smallskip

{\em (b)}
Let $G$ be a perfect graph and 
suppose that the toric ring $K[\QQ_G]$ of the stable polytope 
$\QQ_G$ of $G$ is Gorenstein. 
Then there is a special simplex in $\QQ_G$.
\end{Theorem}

\begin{pf}
(a) Let $[n]$ be the vertex set of $G$.
Since $G$ possesses a perfect matching,
it follows that $n = 2m$ is even and
there exist $m$ edges $e_1, \ldots, e_m$ of $G$
with $\rho(e_1) + \ldots + \rho(e_m)
= \eb_1 + \cdots + \eb_n$.
Let $\Sigma$ denote the $(m - 1)$-simplex whose vertices are
$\rho(e_1), \ldots, \rho(e_m)$.
We claim that $\Sigma$ is special in $\PP_G$.

Theorem \ref{Gorenstein}
together with \cite[Theorem 1.7]{OhHinormal}
give the complete information about 
the equations of the supporting hyperplanes 
which define the facets of the edge polytope $\PP_G$.
Let $\HH_i$ denote the hyperplane defined by the
equation $z_i = 0$.
Then $\rho(e_j) \in \HH_i$ if and only if 
$i \not\in e_j$.
% Hence $\HH_i$ contains exactly
% $m - 1$ vertices of $\Sigma$.
Let $\HH_T$ denote the hyperplane defined by the equation 
$\sum_{i \in T} z_i = \sum_{j \in N(G;T)} z_j$,
where $T$ is a stable set of $G$,
which is the supporting hyperplane of a facet of $\PP_G$.
Let $C$ denote the set of those $1 \leq i \leq m$ with 
$\rho(e_i) \bigcap T \neq \emptyset$ 
and $D$ the set of those $1 \leq i \leq m$ with 
$\rho(e_i) \Sect N(G;T) \neq \emptyset$.
In either the case of $|N(G;T)| = |T| + 1$
or the case of $|N(G;T)| -1 = |T| + 1 = m$,
one has
$C \subset D$ with $|C| = |D| - 1$.
Let $i_0 \in D \setminus C$.
Then $\rho(e_j) \in \HH_T$
if and only if $j \neq i_0$.
% Hence $\HH_T$ contains exactly
% $m - 1$ vertices of $\Sigma$.
Thus $\Sigma$ is special in $\PP_G$
as desired.

% \smallskip

(b)
Let $G$ be a perfect graph on $[n]$
and suppose that 
all maximal cliques have the cardinality $q$.
Since $G$ is perfect, 
the chromatic number of $G$ is equal to $q$.  Thus 
there is a map $\varphi : [n] \to [q]$
with the property that $\varphi(i) \neq \varphi(j)$
if $\{i, j \} \in E(G)$.
Let $W'_\ell$ denote the stable set
$\{ i \in [n] \, : \, \varphi(i) = \ell \}$ 
for each $1 \leq \ell \leq q$.
We assume that $\QQ_G$ is not a simplex.
Thus one of the stable sets $W'_1, \ldots, W'_q$ 
contains at least two vertices.
Let, say, $W'_1$ contain at least two vertices
and fix $i_0 \in W'_1$.
Let $W_0 = \{ i_0 \}, W_1 = W'_1 \setminus \{ i_0 \}$
and $W_\ell = W'_\ell$ for $2 \leq \ell \leq q$.
Each of $W_0, W_1, \ldots, W_q$ is a stable set
of $G$ and 
$\sum_{\ell = 0}^{q} \rho(W_\ell) 
= \eb_1 + \cdots + \eb_n$.
Let $\Sigma$ denote the $q$-simplex with $q + 1$
vertices $\rho(W_0), \rho(W_1), \ldots, \rho(W_q)$. 
We claim that $\Sigma$ is special in $\QQ_G$. 
 
Recall that 
the equation of the supporting hyperplanes which defines 
the facets of the stable polytope of $G$ are 
either (i)
$x_i = 0$ for $1 \leq i \leq n$
or (ii)
$\sum_{i \in W} x_i = 1$, where $W$ is a maximal clique of $G$. 
If $F_i$ is the facet defined by $x_i = 0$, then
$\rho(W_\ell) \in F_i$ if and only if $i \not\in W_\ell$.
Since $[n]$ is the disjoint union
$W_0 \Union W_1 \Union \cdots \Union W_q$,
it follows that $F_i$ contains exactly $q$ of the vertices
of $\Sigma$.
Let $F'_W$ denote the facet defined by 
$\sum_{i \in W} x_i = 1$, where $W$ is a maximal clique 
of $G$.
Since each of the subsets
$W \cap (W_0 \Union W_1), 
W \cap W_2, \ldots, W \cap W_q$
of $[n]$ consists of one element,
it follows that
each of the vertices
$\rho(W_2), \ldots, \rho(W_q)$
belongs to
$F'_W$ and that
$\rho(W_0) \in F'_W$ 
(resp. $\rho(W_1) \in F'_W$)
if and only if
$i_0 \in W$
(resp. $i_0 \not\in W$).
Hence $F'_W$ contains exactly $q$ of the vertices
of $\Sigma$.
\end{pf}

% \newpage

By using Example \ref{posetpolytope} 
together with \cite[Theorem 3.2]{twoposetpolytopes}, 
it turns out that the above Theorem \ref{specialsimplex} (b)
generalize Reiner--Welker \cite[Corollary 3.8]{RW}.

Now, by virtue of \cite[Theorem 3.5]{Christos},
one has 
a rich class of unimodal Gorenstein sequences
\cite[p. 66]{StanleyGreenBook}.
It is known \cite[Example 1.3 (c)]{OhHicompressed}
the stable polytope of a perfect graph is compressed.

The edge polytope of a finite connected graph is unimodular,
i.e., all of its triangulations are unimodular,
if and only if any two odd cycles of $G$ possess at least
one common vertex.  In particular 
the edge polytope of a finite connected bipartite graph is 
unimodular.
The edge polytope of $G_1$ of Example \ref{graph} 
is compressed (\cite{OhHimultipartite}) but not unimodular,
and that of $G_2$
is unimodular.
%It seems, of course, be of interest (however, be of difficult)
%to find 
A combinatorial characterization of finite graphs $G$
for which the edge polytope $\PP_G$ is compressed is
given in \cite[Theorem 4.1]{Ohsugi}.

\begin{Corollary}
{\em (a)}
Let $G$ be a finite connected graph 
as in Theorem \ref{Gorenstein} (a) 
and suppose that the edge polytope $\PP_G$ is compressed 
and that
the toric ring $K[\PP_G]$ is Gorenstein.   
Then the $h$-vector of $K[\PP_G]$ is unimodal.

{\em (a')}
Let $G$ be a finite $2$-connected bipartite graph 
and suppose that the toric ring $K[\PP_G]$ is Gorenstein.   
Then the $h$-vector of $K[\PP_G]$ is unimodal.

% \smallskip

{\em (b)}
Let $G$ be a perfect graph and 
suppose that the toric ring $K[\QQ_G]$ of the stable polytope 
$\QQ_G$ of $G$ is Gorenstein. 
Then the $h$-vector of $K[\QQ_G]$ is unimodal.
\end{Corollary}

We conclude the present paper with

\begin{Example}
\label{final}
{\em
Let $n \geq 3$ and $G$ the finite connected graph 
on $[2n]$ drawn below.
Let $n$ be odd.
(If $n$ is even, then $K[\PP_G]$
is not Gorenstein by Theorem \ref{Gorenstein} (a').)
By virtue of \cite[Theorem 4.1]{Ohsugi}
it turns out that the edge polytope $\PP_G$ of $G$ is 
compressed.  By using Theorem \ref{Gorenstein} (a)
it follows that the toric ring $K[\PP_G]$
is (normal and) Gorenstein.
Moreover, we can compute the $h$-vector explicitly.
Since the graph $G$ satisfies the condition in
\cite[Theorem 1.2]{OhHiquadratic},
the ``toric ideal" $I_G$ of $G$ is generated by quadratic binomials
which correspond to even cycles of $G$
of length $4$.
There exists a reverse lexicographic order
such that the initial monomials of the quadratic binomials 
are relatively prime.
%Let $<_{rev}$ be a reverse lexicographic order
%such that the initial monomial of quadratic binomials is relatively prime.
Since the set of quadratic binomials is a Gr\"obner basis 
of $I_G$ with respect to $<_{rev}$, it follows that
the initial ideal is generated by $n$ monomials 
which are squarefree, quadratic and relatively prime.
Thus the $h$-vector of 
$K[\PP_G]$ is 
$\left(1, n, {n \choose 2}, \cdots, {n \choose n - 2}, n, 1\right)$.
}
\end{Example}

\bigskip

\begin{center}
%WinTpicVersion2.15
\unitlength 0.1in
\begin{picture}(33.35,34.85)(18.65,-36.00)
% CIRCLE 2 0 0 0
% 4 2800 1000 2870 1070 2860 1070 2860 1070
% 
\special{pn 8}%
\special{sh 0.600}%
\special{ar 2800 600 99 99  0.0000000 6.2831853}%
% STR 2 0 3 0
% 3 3600 500 3600 600 5 0
% $1$
\put(36.0000,-2.0000){\makebox(0,0){$1$}}%
% STR 2 0 3 0
% 3 2600 700 2600 800 5 0
% $2$
\put(26.0000,-4.0000){\makebox(0,0){$2$}}%
% STR 2 0 3 0
% 3 2000 1300 2000 1400 5 0
% $3$
\put(20.0000,-10.0000){\makebox(0,0){$3$}}%
% STR 2 0 3 0
% 3 4600 700 4600 800 5 0
% $n$
\put(46.0000,-4.0000){\makebox(0,0){$n$}}%
% CIRCLE 2 0 0 0
% 4 4400 1000 4470 1070 4460 1070 4460 1070
% 
\special{pn 8}%
\special{sh 0.600}%
\special{ar 4400 600 99 99  0.0000000 6.2831853}%
% CIRCLE 2 0 0 0
% 4 3600 800 3670 870 3660 870 3660 870
% 
\special{pn 8}%
\special{sh 0.600}%
\special{ar 3600 400 99 99  0.0000000 6.2831853}%
% CIRCLE 2 0 0 0
% 4 2200 1600 2270 1670 2260 1670 2260 1670
% 
\special{pn 8}%
\special{sh 0.600}%
\special{ar 2200 1200 99 99  0.0000000 6.2831853}%
% CIRCLE 2 0 0 0
% 4 3600 1400 3670 1470 3660 1470 3660 1470
% 
\special{pn 8}%
\special{sh 0.600}%
\special{ar 3600 1000 99 99  0.0000000 6.2831853}%
% CIRCLE 2 0 0 0
% 4 3090 1560 3160 1630 3150 1630 3150 1630
% 
\special{pn 8}%
\special{sh 0.600}%
\special{ar 3090 1160 99 99  0.0000000 6.2831853}%
% CIRCLE 2 0 0 0
% 4 4100 1540 4170 1610 4160 1610 4160 1610
% 
\special{pn 8}%
\special{sh 0.600}%
\special{ar 4100 1140 99 99  0.0000000 6.2831853}%
% CIRCLE 2 0 0 0
% 4 2690 1990 2760 2060 2750 2060 2750 2060
% 
\special{pn 8}%
\special{sh 0.600}%
\special{ar 2690 1590 99 99  0.0000000 6.2831853}%
% CIRCLE 2 2 3 0
% 4 3600 2400 3600 4000 2000 2200 5000 1400
% 
\special{pn 8}%
\special{ar 3600 2000 1600 1600  5.6629358 5.6704358}%
\special{ar 3600 2000 1600 1600  5.6929358 5.7004358}%
\special{ar 3600 2000 1600 1600  5.7229358 5.7304358}%
\special{ar 3600 2000 1600 1600  5.7529358 5.7604358}%
\special{ar 3600 2000 1600 1600  5.7829358 5.7904358}%
\special{ar 3600 2000 1600 1600  5.8129358 5.8204358}%
\special{ar 3600 2000 1600 1600  5.8429358 5.8504358}%
\special{ar 3600 2000 1600 1600  5.8729358 5.8804358}%
\special{ar 3600 2000 1600 1600  5.9029358 5.9104358}%
\special{ar 3600 2000 1600 1600  5.9329358 5.9404358}%
\special{ar 3600 2000 1600 1600  5.9629358 5.9704358}%
\special{ar 3600 2000 1600 1600  5.9929358 6.0004358}%
\special{ar 3600 2000 1600 1600  6.0229358 6.0304358}%
\special{ar 3600 2000 1600 1600  6.0529358 6.0604358}%
\special{ar 3600 2000 1600 1600  6.0829358 6.0904358}%
\special{ar 3600 2000 1600 1600  6.1129358 6.1204358}%
\special{ar 3600 2000 1600 1600  6.1429358 6.1504358}%
\special{ar 3600 2000 1600 1600  6.1729358 6.1804358}%
\special{ar 3600 2000 1600 1600  6.2029358 6.2104358}%
\special{ar 3600 2000 1600 1600  6.2329358 6.2404358}%
\special{ar 3600 2000 1600 1600  6.2629358 6.2704358}%
\special{ar 3600 2000 1600 1600  6.2929358 6.3004358}%
\special{ar 3600 2000 1600 1600  6.3229358 6.3304358}%
\special{ar 3600 2000 1600 1600  6.3529358 6.3604358}%
\special{ar 3600 2000 1600 1600  6.3829358 6.3904358}%
\special{ar 3600 2000 1600 1600  6.4129358 6.4204358}%
\special{ar 3600 2000 1600 1600  6.4429358 6.4504358}%
\special{ar 3600 2000 1600 1600  6.4729358 6.4804358}%
\special{ar 3600 2000 1600 1600  6.5029358 6.5104358}%
\special{ar 3600 2000 1600 1600  6.5329358 6.5404358}%
\special{ar 3600 2000 1600 1600  6.5629358 6.5704358}%
\special{ar 3600 2000 1600 1600  6.5929358 6.6004358}%
\special{ar 3600 2000 1600 1600  6.6229358 6.6304358}%
\special{ar 3600 2000 1600 1600  6.6529358 6.6604358}%
\special{ar 3600 2000 1600 1600  6.6829358 6.6904358}%
\special{ar 3600 2000 1600 1600  6.7129358 6.7204358}%
\special{ar 3600 2000 1600 1600  6.7429358 6.7504358}%
\special{ar 3600 2000 1600 1600  6.7729358 6.7804358}%
\special{ar 3600 2000 1600 1600  6.8029358 6.8104358}%
\special{ar 3600 2000 1600 1600  6.8329358 6.8404358}%
\special{ar 3600 2000 1600 1600  6.8629358 6.8704358}%
\special{ar 3600 2000 1600 1600  6.8929358 6.9004358}%
\special{ar 3600 2000 1600 1600  6.9229358 6.9304358}%
\special{ar 3600 2000 1600 1600  6.9529358 6.9604358}%
\special{ar 3600 2000 1600 1600  6.9829358 6.9904358}%
\special{ar 3600 2000 1600 1600  7.0129358 7.0204358}%
\special{ar 3600 2000 1600 1600  7.0429358 7.0504358}%
\special{ar 3600 2000 1600 1600  7.0729358 7.0804358}%
\special{ar 3600 2000 1600 1600  7.1029358 7.1104358}%
\special{ar 3600 2000 1600 1600  7.1329358 7.1404358}%
\special{ar 3600 2000 1600 1600  7.1629358 7.1704358}%
\special{ar 3600 2000 1600 1600  7.1929358 7.2004358}%
\special{ar 3600 2000 1600 1600  7.2229358 7.2304358}%
\special{ar 3600 2000 1600 1600  7.2529358 7.2604358}%
\special{ar 3600 2000 1600 1600  7.2829358 7.2904358}%
\special{ar 3600 2000 1600 1600  7.3129358 7.3204358}%
\special{ar 3600 2000 1600 1600  7.3429358 7.3504358}%
\special{ar 3600 2000 1600 1600  7.3729358 7.3804358}%
\special{ar 3600 2000 1600 1600  7.4029358 7.4104358}%
\special{ar 3600 2000 1600 1600  7.4329358 7.4404358}%
\special{ar 3600 2000 1600 1600  7.4629358 7.4704358}%
\special{ar 3600 2000 1600 1600  7.4929358 7.5004358}%
\special{ar 3600 2000 1600 1600  7.5229358 7.5304358}%
\special{ar 3600 2000 1600 1600  7.5529358 7.5604358}%
\special{ar 3600 2000 1600 1600  7.5829358 7.5904358}%
\special{ar 3600 2000 1600 1600  7.6129358 7.6204358}%
\special{ar 3600 2000 1600 1600  7.6429358 7.6504358}%
\special{ar 3600 2000 1600 1600  7.6729358 7.6804358}%
\special{ar 3600 2000 1600 1600  7.7029358 7.7104358}%
\special{ar 3600 2000 1600 1600  7.7329358 7.7404358}%
\special{ar 3600 2000 1600 1600  7.7629358 7.7704358}%
\special{ar 3600 2000 1600 1600  7.7929358 7.8004358}%
\special{ar 3600 2000 1600 1600  7.8229358 7.8304358}%
\special{ar 3600 2000 1600 1600  7.8529358 7.8604358}%
\special{ar 3600 2000 1600 1600  7.8829358 7.8904358}%
\special{ar 3600 2000 1600 1600  7.9129358 7.9204358}%
\special{ar 3600 2000 1600 1600  7.9429358 7.9504358}%
\special{ar 3600 2000 1600 1600  7.9729358 7.9804358}%
\special{ar 3600 2000 1600 1600  8.0029358 8.0104358}%
\special{ar 3600 2000 1600 1600  8.0329358 8.0404358}%
\special{ar 3600 2000 1600 1600  8.0629358 8.0704358}%
\special{ar 3600 2000 1600 1600  8.0929358 8.1004358}%
\special{ar 3600 2000 1600 1600  8.1229358 8.1304358}%
\special{ar 3600 2000 1600 1600  8.1529358 8.1604358}%
\special{ar 3600 2000 1600 1600  8.1829358 8.1904358}%
\special{ar 3600 2000 1600 1600  8.2129358 8.2204358}%
\special{ar 3600 2000 1600 1600  8.2429358 8.2504358}%
\special{ar 3600 2000 1600 1600  8.2729358 8.2804358}%
\special{ar 3600 2000 1600 1600  8.3029358 8.3104358}%
\special{ar 3600 2000 1600 1600  8.3329358 8.3404358}%
\special{ar 3600 2000 1600 1600  8.3629358 8.3704358}%
\special{ar 3600 2000 1600 1600  8.3929358 8.4004358}%
\special{ar 3600 2000 1600 1600  8.4229358 8.4304358}%
\special{ar 3600 2000 1600 1600  8.4529358 8.4604358}%
\special{ar 3600 2000 1600 1600  8.4829358 8.4904358}%
\special{ar 3600 2000 1600 1600  8.5129358 8.5204358}%
\special{ar 3600 2000 1600 1600  8.5429358 8.5504358}%
\special{ar 3600 2000 1600 1600  8.5729358 8.5804358}%
\special{ar 3600 2000 1600 1600  8.6029358 8.6104358}%
\special{ar 3600 2000 1600 1600  8.6329358 8.6404358}%
\special{ar 3600 2000 1600 1600  8.6629358 8.6704358}%
\special{ar 3600 2000 1600 1600  8.6929358 8.7004358}%
\special{ar 3600 2000 1600 1600  8.7229358 8.7304358}%
\special{ar 3600 2000 1600 1600  8.7529358 8.7604358}%
\special{ar 3600 2000 1600 1600  8.7829358 8.7904358}%
\special{ar 3600 2000 1600 1600  8.8129358 8.8204358}%
\special{ar 3600 2000 1600 1600  8.8429358 8.8504358}%
\special{ar 3600 2000 1600 1600  8.8729358 8.8804358}%
\special{ar 3600 2000 1600 1600  8.9029358 8.9104358}%
\special{ar 3600 2000 1600 1600  8.9329358 8.9404358}%
\special{ar 3600 2000 1600 1600  8.9629358 8.9704358}%
\special{ar 3600 2000 1600 1600  8.9929358 9.0004358}%
\special{ar 3600 2000 1600 1600  9.0229358 9.0304358}%
\special{ar 3600 2000 1600 1600  9.0529358 9.0604358}%
\special{ar 3600 2000 1600 1600  9.0829358 9.0904358}%
\special{ar 3600 2000 1600 1600  9.1129358 9.1204358}%
\special{ar 3600 2000 1600 1600  9.1429358 9.1504358}%
\special{ar 3600 2000 1600 1600  9.1729358 9.1804358}%
\special{ar 3600 2000 1600 1600  9.2029358 9.2104358}%
\special{ar 3600 2000 1600 1600  9.2329358 9.2404358}%
\special{ar 3600 2000 1600 1600  9.2629358 9.2704358}%
\special{ar 3600 2000 1600 1600  9.2929358 9.3004358}%
\special{ar 3600 2000 1600 1600  9.3229358 9.3304358}%
\special{ar 3600 2000 1600 1600  9.3529358 9.3604358}%
\special{ar 3600 2000 1600 1600  9.3829358 9.3904358}%
\special{ar 3600 2000 1600 1600  9.4129358 9.4204358}%
\special{ar 3600 2000 1600 1600  9.4429358 9.4504358}%
\special{ar 3600 2000 1600 1600  9.4729358 9.4804358}%
\special{ar 3600 2000 1600 1600  9.5029358 9.5104358}%
\special{ar 3600 2000 1600 1600  9.5329358 9.5404358}%
% CIRCLE 2 0 3 0
% 4 3600 2400 2000 2200 2140 1690 2000 2200
% 
\special{pn 8}%
\special{ar 3600 2000 1612 1612  3.2659476 3.5942214}%
% CIRCLE 2 0 3 0
% 4 3600 2400 3600 800 2720 1060 2270 1550
% 
\special{pn 8}%
\special{ar 3600 2000 1600 1600  3.7102655 4.1313067}%
% CIRCLE 2 0 3 0
% 4 3600 2400 3600 800 3500 800 2890 950
% 
\special{pn 8}%
\special{ar 3600 2000 1600 1600  4.2570514 4.6499702}%
% CIRCLE 2 0 3 0
% 4 3600 2400 3600 800 4310 950 3700 780
% 
\special{pn 8}%
\special{ar 3600 2000 1600 1600  4.7740392 5.1677265}%
% CIRCLE 2 0 3 0
% 4 3600 2400 3600 800 5000 1600 4480 1050
% 
\special{pn 8}%
\special{ar 3600 2000 1600 1600  5.2900649 5.7640392}%
% CIRCLE 2 2 3 0
% 4 3600 2400 3600 1400 2600 2400 4400 1800
% 
\special{pn 8}%
\special{ar 3600 2000 1000 1000  5.6396842 5.6516842}%
\special{ar 3600 2000 1000 1000  5.6876842 5.6996842}%
\special{ar 3600 2000 1000 1000  5.7356842 5.7476842}%
\special{ar 3600 2000 1000 1000  5.7836842 5.7956842}%
\special{ar 3600 2000 1000 1000  5.8316842 5.8436842}%
\special{ar 3600 2000 1000 1000  5.8796842 5.8916842}%
\special{ar 3600 2000 1000 1000  5.9276842 5.9396842}%
\special{ar 3600 2000 1000 1000  5.9756842 5.9876842}%
\special{ar 3600 2000 1000 1000  6.0236842 6.0356842}%
\special{ar 3600 2000 1000 1000  6.0716842 6.0836842}%
\special{ar 3600 2000 1000 1000  6.1196842 6.1316842}%
\special{ar 3600 2000 1000 1000  6.1676842 6.1796842}%
\special{ar 3600 2000 1000 1000  6.2156842 6.2276842}%
\special{ar 3600 2000 1000 1000  6.2636842 6.2756842}%
\special{ar 3600 2000 1000 1000  6.3116842 6.3236842}%
\special{ar 3600 2000 1000 1000  6.3596842 6.3716842}%
\special{ar 3600 2000 1000 1000  6.4076842 6.4196842}%
\special{ar 3600 2000 1000 1000  6.4556842 6.4676842}%
\special{ar 3600 2000 1000 1000  6.5036842 6.5156842}%
\special{ar 3600 2000 1000 1000  6.5516842 6.5636842}%
\special{ar 3600 2000 1000 1000  6.5996842 6.6116842}%
\special{ar 3600 2000 1000 1000  6.6476842 6.6596842}%
\special{ar 3600 2000 1000 1000  6.6956842 6.7076842}%
\special{ar 3600 2000 1000 1000  6.7436842 6.7556842}%
\special{ar 3600 2000 1000 1000  6.7916842 6.8036842}%
\special{ar 3600 2000 1000 1000  6.8396842 6.8516842}%
\special{ar 3600 2000 1000 1000  6.8876842 6.8996842}%
\special{ar 3600 2000 1000 1000  6.9356842 6.9476842}%
\special{ar 3600 2000 1000 1000  6.9836842 6.9956842}%
\special{ar 3600 2000 1000 1000  7.0316842 7.0436842}%
\special{ar 3600 2000 1000 1000  7.0796842 7.0916842}%
\special{ar 3600 2000 1000 1000  7.1276842 7.1396842}%
\special{ar 3600 2000 1000 1000  7.1756842 7.1876842}%
\special{ar 3600 2000 1000 1000  7.2236842 7.2356842}%
\special{ar 3600 2000 1000 1000  7.2716842 7.2836842}%
\special{ar 3600 2000 1000 1000  7.3196842 7.3316842}%
\special{ar 3600 2000 1000 1000  7.3676842 7.3796842}%
\special{ar 3600 2000 1000 1000  7.4156842 7.4276842}%
\special{ar 3600 2000 1000 1000  7.4636842 7.4756842}%
\special{ar 3600 2000 1000 1000  7.5116842 7.5236842}%
\special{ar 3600 2000 1000 1000  7.5596842 7.5716842}%
\special{ar 3600 2000 1000 1000  7.6076842 7.6196842}%
\special{ar 3600 2000 1000 1000  7.6556842 7.6676842}%
\special{ar 3600 2000 1000 1000  7.7036842 7.7156842}%
\special{ar 3600 2000 1000 1000  7.7516842 7.7636842}%
\special{ar 3600 2000 1000 1000  7.7996842 7.8116842}%
\special{ar 3600 2000 1000 1000  7.8476842 7.8596842}%
\special{ar 3600 2000 1000 1000  7.8956842 7.9076842}%
\special{ar 3600 2000 1000 1000  7.9436842 7.9556842}%
\special{ar 3600 2000 1000 1000  7.9916842 8.0036842}%
\special{ar 3600 2000 1000 1000  8.0396842 8.0516842}%
\special{ar 3600 2000 1000 1000  8.0876842 8.0996842}%
\special{ar 3600 2000 1000 1000  8.1356842 8.1476842}%
\special{ar 3600 2000 1000 1000  8.1836842 8.1956842}%
\special{ar 3600 2000 1000 1000  8.2316842 8.2436842}%
\special{ar 3600 2000 1000 1000  8.2796842 8.2916842}%
\special{ar 3600 2000 1000 1000  8.3276842 8.3396842}%
\special{ar 3600 2000 1000 1000  8.3756842 8.3876842}%
\special{ar 3600 2000 1000 1000  8.4236842 8.4356842}%
\special{ar 3600 2000 1000 1000  8.4716842 8.4836842}%
\special{ar 3600 2000 1000 1000  8.5196842 8.5316842}%
\special{ar 3600 2000 1000 1000  8.5676842 8.5796842}%
\special{ar 3600 2000 1000 1000  8.6156842 8.6276842}%
\special{ar 3600 2000 1000 1000  8.6636842 8.6756842}%
\special{ar 3600 2000 1000 1000  8.7116842 8.7236842}%
\special{ar 3600 2000 1000 1000  8.7596842 8.7716842}%
\special{ar 3600 2000 1000 1000  8.8076842 8.8196842}%
\special{ar 3600 2000 1000 1000  8.8556842 8.8676842}%
\special{ar 3600 2000 1000 1000  8.9036842 8.9156842}%
\special{ar 3600 2000 1000 1000  8.9516842 8.9636842}%
\special{ar 3600 2000 1000 1000  8.9996842 9.0116842}%
\special{ar 3600 2000 1000 1000  9.0476842 9.0596842}%
\special{ar 3600 2000 1000 1000  9.0956842 9.1076842}%
\special{ar 3600 2000 1000 1000  9.1436842 9.1556842}%
\special{ar 3600 2000 1000 1000  9.1916842 9.2036842}%
\special{ar 3600 2000 1000 1000  9.2396842 9.2516842}%
\special{ar 3600 2000 1000 1000  9.2876842 9.2996842}%
\special{ar 3600 2000 1000 1000  9.3356842 9.3476842}%
\special{ar 3600 2000 1000 1000  9.3836842 9.3956842}%
% CIRCLE 2 0 3 0
% 4 3600 2400 3600 1400 2640 2070 2600 2400
% 
\special{pn 8}%
\special{ar 3600 2000 1000 1000  3.1415927 3.4726887}%
% CIRCLE 2 0 3 0
% 4 3600 2400 3600 1400 3490 1400 3170 1480
% 
\special{pn 8}%
\special{ar 3600 2000 1000 1000  4.2751669 4.6028295}%
% CIRCLE 2 0 3 0
% 4 3600 2400 3600 1400 2990 1600 2740 1910
% 
\special{pn 8}%
\special{ar 3600 2000 1000 1000  3.6594856 4.0609358}%
% CIRCLE 2 0 3 0
% 4 3600 2400 3600 1400 4020 1490 3690 1410
% 
\special{pn 8}%
\special{ar 3600 2000 1000 1000  4.8030489 5.1447968}%
% CIRCLE 2 0 3 0
% 4 3600 2400 3600 1400 4400 1800 4190 1600
% 
\special{pn 8}%
\special{ar 3600 2000 1000 1000  5.3478420 5.6396842}%
% STR 2 0 3 0
% 3 3600 1500 3600 1600 5 0
% $n+1$
\put(36.0000,-12.0000){\makebox(0,0){$n+1$}}%
% STR 2 0 3 0
% 3 3210 1730 3210 1830 5 0
% $n+2$
\put(32.1000,-14.3000){\makebox(0,0){$n+2$}}%
% STR 2 0 3 0
% 3 2950 2080 2950 2180 5 0
% $n+3$
\put(29.5000,-17.8000){\makebox(0,0){$n+3$}}%
% STR 2 0 3 0
% 3 4000 1700 4000 1800 5 0
% $2 n$
\put(40.0000,-14.0000){\makebox(0,0){$2 n$}}%
% LINE 2 0 3 0
% 2 3600 890 3600 1300
% 
\special{pn 8}%
\special{pa 3600 490}%
\special{pa 3600 900}%
\special{fp}%
% LINE 2 0 3 0
% 2 2840 1080 3030 1460
% 
\special{pn 8}%
\special{pa 2840 680}%
\special{pa 3030 1060}%
\special{fp}%
% LINE 2 0 3 0
% 2 2270 1650 2610 1930
% 
\special{pn 8}%
\special{pa 2270 1250}%
\special{pa 2610 1530}%
\special{fp}%
% LINE 2 0 3 0
% 2 4370 1090 4160 1470
% 
\special{pn 8}%
\special{pa 4370 690}%
\special{pa 4160 1070}%
\special{fp}%
\end{picture}%
\end{center}

\bigskip

{\small
\noindent
Hidefumi Ohsugi
\hfill 
Takayuki Hibi\\
Department of Mathematics
\hfill
Department of Pure and Applied Mathematics\\
Faculty of Science
\hfill
Graduate School of Information Science and Technology\\
Rikkyo University
\hfill
Osaka University \\
Toshima, Tokyo 171--8501, Japan
\hfill
Toyonaka, Osaka 560--0043, Japan\\
E-mail:ohsugi@@rkmath.rikkyo.ac.jp
\hfill
E-mail:hibi@@math.sci.osaka-u.ac.jp

}

\end{document}